\documentclass[12pt]{article}

\usepackage{epsfig}
\usepackage{latexsym}
\usepackage{pdfsync}

\newcommand{\eq}{\begin{equation}\begin{array}{rllllllllllllllllllllllllllllllll}}
\newcommand{\ee}{\end{array}\end{equation}}
\newcommand{\bmt}{\left[ \begin{array}{ccccccccc}}
\newcommand{\emt}{\end{array}\right]}
\newcommand{\bea}{\begin{eqnarray}}

\newcommand{\eea}{\end{eqnarray}}
\newcommand{\bean}{\begin{eqnarray*}}
\newcommand{\eean}{\end{eqnarray*}}

\title{Series Solution of Discrete Time \\Stochastic Optimal Control Problems}
\author{Arthur J Krener
\thanks{Research supported in part by AFOSR .}
\thanks{A. J. Krener is with the Department of Applied Mathematics, Naval Postgraduate School, Monterey, CA 93943
        {\tt\small ajkrener@nps.edu}}}%

\begin{document}

\date{}

\maketitle

\begin{abstract}
 In this paper we consider discrete time stochastic optimal control problems over infinite and finite time horizons.  We show that for a large class of such problems the Taylor polynomials of  the solutions to the associated Dynamic Programming Equations can be  computed degree by degree.
  \end{abstract}

\section{ Introduction}
\setcounter{equation}{0}
We begin with a relatively simple stochastic infinite horizon optimal control problem and then move on to more complicated problems over infinite and finite horizons.
Consider a discrete time, infinite horizon,  stochastic Linear Quadratic Regulator with Bilinear Noise (DLQGB),
\bean
\min_{u(\cdot)} {1\over 2}{\rm E}\left\{ \sum_0^\infty \ \left(x'Qx+2x'Su+u'Ru\right) \right\}
\eean
subject to
\bean
x^+&=& Fx+Gu+ \sum_{k=1}^r w_k(C_{k} x+D_k u )\\
x(0)&=&x^0
\eean
where $x^+(t)=x(t+1)$.


The state $x$ is $n$ dimensional, the control $u$ is $m$ dimensional and $w(t)=(w_1(t),\ldots, w_r(t))'$ is $r$ dimensional sequence of independent Gaussian random vectors of mean zero and covariance $I^{r\times r}$.
The matrices are sized accordingly, in particular $C_{k}$ is an $n\times n$ matrix and $D_k$ is an $n \times m$ matrix  for each $k=1,\ldots,r$.

To the best of our knowledge discrete time infinite horizon  problems with bilinear noise have not been considered before. In \cite{Kr18} we studied the continuous time version of this problem.   The finite horizon version of this problem with noise entering linearly  is well studied in both discrete \cite{Be17}
and continuous time \cite{FR75}, \cite{YZ99}.  

 We restrict our attention to problems with bilinear noise so that we can use power series techniques to solve the dynamic programming equations of   nonlinear optimal control problems.  The class of infinite horizon nonlinear optimal control problems that  are of interest are of the form
 \bean
\min_{u(\cdot)} {\rm E}\left\{\sum_0^\infty  l(x,u)\right\}
\eean
subject to
\bean
x^+&=& f(x,u)+ \sum_{k=1}^r w_k\gamma_k(x,u)\\
x(0)&=&x^0
\eean
 where $x^+(t)=x(t+1)$, $f(x,u)$ and $\gamma_k(x,u)$ are smooth functions of order $O(x,u)$ and $l(x,u)$ is a smooth function of order $O(x,u)^2$.
 
 Associated to these problems are dynamic programming equations for the optimal cost and optimal feedback.   Assuming they exist, let  $\pi(x)$ be the optimal cost starting at $x$ and  $u=\kappa(x)$ be the optimal feedback at $x$ for this problem.  Then they satisfy the Stochastic  Infinite Horizon Dynamic Programming Equations (SIDPE),
\bea \label{SIDPE1}
\pi(x)&=&  \mbox{min}_u {\rm E} \left\{ \pi(f(x,u)+ \sum_{k=1}^r w_k\gamma_k(x,u))+l(x,u) \right\}\\
\label{SIDPE2}
\kappa(x)&=& \mbox{argmin}_u {\rm E} \left\{ \pi(f(x,u)+ \sum_{k=1}^r w_k\gamma_k(x,u))+l(x,u) \right\}
\eea
These equations differ from their deterministic counteparts because of the presence of the noise terms.

 The class of finite horizon nonlinear optimal control problems that  are of interest are of the form
\bean
\min_{u(\cdot)} {\rm E}\left\{\sum_0^T  l(t,x,u) +\pi_T(x(T))\right\}
\eean
subject to
\bean
x^+&=& f(t, x,u)+ \sum_{k=1}^r w_k\gamma_k(t,x,u)\\
x(t_0)&=&x^0
\eean
 where $f(t,x,u)$ and $\gamma_k(t, x,u)$ are smooth $n$ vector vaued functions with respect to $x,u$ of order $O(x,u)$ and continuous with respect to $t$,  $l(t,x,u)$ is a smooth scalar valued function with respect to $x,u$ of order $O(x,u)^2$ and continuous with respect to $t$ and $\pi_T(x)$ is a smooth function with respect to $x$ of order $O(x)^2$.

 Assuming they exist, let  $\pi(t_0,x^0)$ be the optimal cost given that $x(t_0)=x^0$ and  $u(t)=\kappa(t,x)$ be the optimal feedback for this problem.
 Then they satisfy the Stochastic Finite Horizon  Dynamic Programming  Equations (SFDPE),
\bea \nonumber
\pi(t_0,x^0)&=& \mbox{min}_{u(t_0)} {\rm E} \left\{ \pi\left(t_0+1,z^0\right) +l(t_0,x^0,u(t_0)) \right\}\label{sdpe1} \\
\nonumber
\kappa(t_0,x^0)&=& \mbox{argmin}_{u(t_0)} {\rm E} \left\{ \pi\left(t_0+1,z^0\right) +l(t_0,x^0,u(t_0)) \right\}\label{sdpe2} 
\eea
where $z^0$ is the random vector
\bean
z^0&=&f(t_0,x^0,u(t_0))+ \sum_{k=1}^r w_k(t_0)\gamma_k(t_0,x^0,u(t_0))
\eean
Again these equations differ from their deterministic counteparts because  of the noise terms.

 The rest of the this paper is organized as follows.   
 In the next section we solve infinite horizon discrete time linear quadratic regulator problems with bilinear noise (DLQGB).   In this case the SIDPE reduces to   stochastic discrete time algebraic Riccati  equations (SDARE). To our knowledge these SDARE are new.   We present an iterative method 
for solving SDARE using a solver for the corresponding deterministic algebraic Riccati equation (DARE) such as MATLAB's dare.m.   This iteration  may or may not converge depending on the relative size of the noise coefficients.  In Section 3 we show how the Taylor polynomials of 
the optimal cost $\pi(x)$  and the optimal feedback $u=\kappa(x)$ of the solution of (SIDPE) (\ref{SIDPE1}, \ref{SIDPE2}) can be computed degree by degree up to the degree of smoothness of the problem.

\section{ Discrete Time Linear Quadratic Regulator\\ with Bilinear Noise}

If we can find a smooth scalar valued function $\pi(x)$ and a smooth $m$ vector valued $\kappa(x)$ satisfying the  Infinite Horizon Stochastic Dynamic Programming Equations (SIDPE) (\ref{SIDPE1}, \ref{SIDPE2}) then by a standard verification argument \cite{FR75} one can show that 
 $\pi(x^0)$ is the optimal cost of starting at $x^0$ and $u(0)=\kappa(x^0)$ is the optimal control at $x^0$.

We make the standard assumptions of deterministic LQR,
\begin{enumerate}
\item The matrix
\bean
\bmt Q&S\\S'&R\emt
\eean
is nonnegative definite.
\item The matrix $R$ 
is positive definite.
\item The pair $F$, $G$ is stabilizable.
\item The pair $Q^{1/2}$, $F$ is detectable where $Q=(Q^{1/2})'Q^{1/2}$.
\end{enumerate}

Because of the linear dynamics and quadratic cost, 
we expect that $\pi(x) $ is a quadratic function of $x$ and $\kappa(x)$ is a linear function of $x$,
\bean
\pi(x)= {1\over 2}x'Px,&&
\kappa(x)= Kx
\eean
We plug these expressions  into SIDPE  and they simplify to
\bea
P&=&  Q+K'RK+(F+GK)'P(F+GK)   \label{sdare1}\\
&&+\sum_{k=1^r}(C_k+D_kK)'P(C_k+D_kK) \nonumber
\\
K&=&-\left(R+G'PG+\sum_{k=1}^rD'_kPD_k\right)^{-1}\label{sdare2}\\
&& \times \left(G'PF+S'+\sum_{k=1}^rD'_kPC_k\right) 
\nonumber
\eea

We call these equations  (\ref{sdare1}, \ref{sdare2}) the Stochastic Discrete Time  Algebraic Riccati Equations (SDARE).
They reduce to the deterministic  Discrete Time  Algebraic Riccati Equations (DARE)  if $C_k=0$ and $D_k=0$ for $k=1,\ldots,r$.

Here is an iterative method for solving SDARE.  Let $P_{(0)}$ be the solution  of the first discrete time deterministic algebraic Riccati equation DARE 
 \bean
0&=&  P_{(0)}F+F'P_{(0)}+Q-(P_{(0)}G+S)R^{-1}(G'P_{(0)}+S') 
\eean
and $K_{(0)}$ be solution of the second  deterministic DARE 
\bean
K_{(0)}&=&-R^{-1}(G'P+S')
\eean

Given $P_{(\tau-1)}$  define
\bean
Q_{(\tau)}&=& Q+\sum_{k=1}^r C'_k P_{(\tau-1)}C_k\\
R_{(\tau)}&=& R+\sum_{k=1}^r D'_k P_{(\tau-1)}D_k\\
S_{(\tau)}&=& S+\sum_{k=1}^r C'_k P_{(\tau-1)}D_k
\eean
Let
 $P_{(\tau)}$ be the solution  of
\bean
0&=& P_{(\tau)}F+F'P_{(\tau)}+Q_{(\tau)}-(P_{(\tau)}G+S_{(\tau)})R_{(\tau)}^{-1}(G'P_{(\tau)}+S'_{(\tau)})
\eean
and
\bean
K_{(\tau)}&=&-R_{(\tau)}^{-1}\left(G'P_{(\tau)}+S_{(\tau)}'\right) 
\eean

If the iteration on $P_{(\tau)}$   converges, that is, for some $\tau$, $P_{(\tau)}\approx P_{(\tau-1)}$ then $P_{(\tau)}$ and $ K_{(\tau)}$ are approximate solutions to SDARE

The solution $P$ of the DARE is the kernel of the optimal cost of a deterministic LQR and since   
\bean
\bmt Q& S\\ S'& R\emt \le \bmt Q_{(\tau-1)}& S_{(\tau-1)}\\S'_{(\tau-1)}& R_{(\tau-1)}\emt \le \bmt Q_{(\tau)}& S_{(\tau)}\\S'_{(\tau)}& R_{(\tau)}\emt
\eean
it follows that $P_{(0)}\le P_{(\tau-1)} \le P_{(\tau)} $,  the iteration is monotonically increasing.
 Computationally we have found  that  if matrices $C_k$ and $D_k$ are not too big relative to $F, G,Q,R,S$ then the iteration conveges.  But if the $C_k$ and $D_k$ are about the same size as $F$ and $G$ or larger the iteration can diverge.    Further study of this issue is needed.   The iteration does converge in the simple example in the next section.

It  is well-known \cite{YZ99} that the first and second standard assumptions of LQR can be violated in a stochastic optimal control problem and still the optimal cost can be finite and positive.    This is true for some SLQRB problems and the reason why can be seen in the above iteration.   For some $\tau^*>0$ 
it may happen that
\bean 
\bmt Q_{(\tau^*)} & S_{(\tau^*)}\\S'_{(\tau^*)}&R_{(\tau^*)}\emt &\ge&0\\
R_{(\tau^*)} &>&0
\eean
then this will happen for all $\tau >\tau^*$
even though this might  not be true when $\tau=0$.   The MATLAB function {\tt dare} does require that the first two LQR assumptions hold so it can be used in the above iteration.

\section{DLQGB Example}
\setcounter{equation}{0}
Here is a simple example with $n=2,m=1,r=2$.
\bean
\min_u {1\over 2}\sum_0^\infty \|x\|^2+u^2\ dt
\eean
subject to 
\bean
x_1^+&=&x_1+ 0.1 x_2 +0.1  w_1 x_1 \\
x_2^+&=&x_2+0.1 u +0.1 w_2(x_2 +u)
\eean
In other words 
\bean
Q=\bmt 1&0\\0&1\emt,& S=\bmt 0\\0\emt, &R=1\\
F= \bmt 1&0.1\\0&1\emt,& G=\bmt 0\\0.1\emt& \\
C_1=\bmt 0.1&0\\0&0\emt, &C_2=\bmt 0&0\\0&0.1\emt&\\
D_1=\bmt 0\\0\emt, &D_2=\bmt 0\\0.1\emt&
\eean

The solution of   the noiseless  DARE is
\bean
P&=& \bmt 18.3422 & 10.9046\\10.9046&18.9110\emt
\\
K&=&-\bmt 0.9170&1.6821\emt
\eean
The eigenvalues of the noiseless  closed loop matrix $F+GK$ are $0.9054\pm0.0443 i $ and are of norm $0.9065$.

The above iteration essentially converges to the solution of the  SDARE in about twenty iterations, the solution is
\bean
P&=& \bmt 22.3884 & 13.2764\\13.2764&21.6311\emt\\
\\
K&=&\bmt -1.3276&-2.1631\emt
\eean
The eigenvalues of the noisy closed loop matrix $F+GK$ are $0.8918 \pm 0.0397i $ and are of norm $0.8927$.

As expected the noisy system is  more difficult to control than the noiseless system and the poles are smaller in norm.  It should be noted that the above iteration diverged to infinity
when the noise coefficients were increased from $0.1$ to $1$.

\section{Nonlinear Stochastic Infinite Horizon DPE}
Suppose the problem is not linear-quadratic,  the dynamics is given by a nonlinear stochastic difference equation
\bean
x^+&=& f(x,u) +\sum_{k=1}^r w_k\gamma_k(x,u)
\eean
 and the criterion to be minimized is
 \bean
\min_{u(\cdot)} {\rm E} \left\{\sum_0^\infty l(x,u) \right\}
\eean
As before the noise  $w(t)=(w_1,\ldots,w_r)'$ is a sequence of independent Gaussian vectors of zero mean and covariance $I^{r\times r}$.

We assume that $f(x,u),\gamma_k(x,u), l(x,u) $ are smooth functions that have Taylor polynomial expansions
around $x=0,u=0$.  We also assume that  $f(x,u)=O(x,u)$, $\gamma_k(x,u)=O(x,u)$ and $l(x,u)=O(x,u)^2$ so
\bean
f(x,u)&=& Fx+Gu+f^{[2]}(x,u)+\ldots+f^{[d]}(x,u)+O(x,u)^{d+1}\\
\gamma_k(x,u)&=& C_k x+D_k u+\gamma_k^{[2]}(x,u)+\ldots+\gamma_k^{[d]}(x,u)+O(x,u)^{d+1}\\
l(x,u)&=&{1\over 2}\left(x'Qx+2x'Su+u'Ru\right) +l^{[3]}(x,u)+\ldots\\&&+l^{[d+1]}(x,u)+O(x,u)^{d+2}
\eean
where $^{[d]}$ indicates the homogeneous polynomial terms of degree $d$.

  Then if they exist  the optimal cost $\pi(x)$ and optimal feedback $u=\kappa(x)$ satisfy SIDPE (\ref{SIDPE1}. \ref{SIDPE2}).  The quantity to be minimized  is a smooth function of $u$ hence  (\ref{SIDPE1}. \ref{SIDPE2}) imply
  \bea \label{SIDPE3}
\pi(x)&=& {\rm E} \left\{ \pi\left(f(x,\kappa(x))+\sum_{k=1}^rw_k\gamma_k(x,\kappa(x))\right)\right\}+  l(x,\kappa(x)) \\
\nonumber
0&=& {\rm E} \left\{\frac{\partial \pi}{\partial x}\left(f(x,\kappa(x))+\sum_{k=1}^rw_k\gamma_k(x,\kappa(x)) \right)\right.
\\&& \left.  \times
\left(\frac{\partial f}{\partial u}(x,\kappa(x)) +\sum_k w_k\frac{\partial  \gamma_k}{\partial u}(x,\kappa(x))\right)\right\}\nonumber\\
&&+\frac{\partial l}{\partial u}(x,\kappa(x)) \label{SIDPE4}
\eea
Of course the reverse implication is not necessarily true  as the quantity to be minimized could have local minima  or stationary points.

   We assume that the optimal cost  and optimal feedback  have similar Taylor polynomial expansions
  \bean
\pi(x)&=& {1\over 2}x'Px +\pi^{[3]}(x)+\ldots+\pi^{[d+1]}(x)+O(x)^{d+2}\\
\kappa(x)&=&Kx+\kappa^{[2]}(x)+\ldots+\kappa^{[d]}(x)+O(x)^{d+1}
\eean
We  plug all these expansions into  equations (\ref{SIDPE3}, \ref{SIDPE4}). 
At lowest degrees, degree two in (\ref{SIDPE3}) and degree one in (\ref{SIDPE4}) we get the familiar SDARE
 (\ref{sdare1}, \ref{sdare2}).

If  (\ref{sdare1}, \ref{sdare2}) are solvable 
then  we may proceed to the next degrees, degree three in (\ref{SIDPE3})  and degree two in (\ref{SIDPE4}).  
\bea 
\label{SIDPE5}
&&\pi^{[3]}(x)={\rm E} \left\{\pi^{[3]}\left( (F+GK)x+\sum_k w_k (C_k+D_kK)x\right)\right\}\\
&&+x'(F+GK)'Pf^{[2]}(x,Kx)+\sum_k x'(C_k+D_kK)'P\gamma_k^{[2]}(x,Kx)+l^{[3]}(x,Kx)\nonumber \\
\label{SIDPE6}
&&0={\rm E} \left\{\frac{\partial \pi^{[3]}}{\partial x}\left( (F+GK)x+\sum_k w_k(C_k+D_kK)x\right)
\left(G+\sum_kw_k D_k  \right)  \right\}\nonumber
\\
\nonumber 
&&+x' (F+GK)'P\frac{\partial f^{[2]}}{\partial u}(x,Kx)+\frac{\partial l^{[3]}}{\partial u}(x,Kx)\\
&&+(\kappa^{[2]}(x))'\left(R+G'PG+\sum_kD_k'PD_k\right)
\eea

Notice the first equation (\ref{SIDPE5})  is a square linear equation for the  unknown $\pi^{[3]}(x)$,
the other unknown $\kappa^{[2]}(x)$ does not appear in it.
 If we can solve  it  for $\pi^{[3]}(x)$
 then we can solve the second equation (\ref{SIDPE5})  for $\kappa^{[2]}(x)$ because of the standard assumption that $R$ is invertible so $R+G'PG+\sum_k D_kPD_k$
 must also be invertible.
 
 In the deterministic case the eigenvalues of the linear operator
 \bea \label{dop}
 \pi^{[3]}(x) &\mapsto&  \pi^{[3]}\left((F+GK)x\right)
 \eea 
 are the products of three eigenvalues of $F+GK$.  Under the standard LQR assumptions all the  eigenvalues of $F+GK$ are in the open unit disc so any
product of three eigenvalues of $F+GK$ has norm less than one.   Hence  the operator   
  \bea \label{dop1}
 \pi^{[3]}(x) &\mapsto& \pi^{[3]}(x)-  \pi^{[3]}\left((F+GK)x\right)
 \eea 
is invertible.
If the noise coefficients $C_k,\ D_k$ are small relative to the eigenvalues of (\ref{dop})  then the operator 
 \bea \label{dop2}
 \pi^{[3]}(x) &\mapsto& \pi^{[3]}(x)-  {\rm E} \left\{\pi^{[3]}\left( (F+GK)x+\sum_k w_k(C_k+D_kK)x\right)\right\}
 \eea 
 will also be invertible and so we can solve (\ref{SIDPE5}) for $\pi^{[3]}(x)$ and then (\ref{SIDPE6}) for $\kappa^{[2]}(x)$.

 The first SIDPE equation  for $\pi^{[d+1]}(x)$ contains previously computed lower degree terms and the linear operator 
  \bea \label{dopd1}
 \pi^{[  d+1]}(x) \mapsto \pi^{[  d+1]}(x)-  {\rm E} \left\{\pi^{[  d+1]}\left( (F+GK)x+\sum_k (C_k+D_kK)x\ w_k\right)\right\}
 \eea 
 The eigenvalues of deterministic part of this operator 
  \bea \label{dopd2}
 \pi^{[  d+1]}(x) &\mapsto& \pi^{[  d+1]}(x)-  \pi^{[  d+1]}\left((F+GK)x\right)
 \eea
 are of the form $1-\lambda_{i_1}\cdots \lambda_{i_{d+1}}$ where $\lambda_{i_j}$  are eigenvalues of $F+Gk$ which are strictly inside the
 unit disk.  Hence (\ref{dopd2}) is always invertible and its stochastic perturbation  (\ref{dopd1}) will be also if $C_k$ and $D_k$ are small enough.

\section{Nonlinear Example}
\setcounter{equation}{0}
Here is a simple example with $n=2,m=1,r=2$.  Consider a pendulum of length $1\ m$ and mass $1\ kg$ orbiting approximately 400 kilometers
above Earth on the International Space Station (ISS).   The "gravity constant" at this height is approximately $g=8.7\ m/sec^2$.  The pendulum can be controlled
by a torque $u$ that can be applied at the pivot and there is damping at the pivot with  linear damping constant $c_1=0.1\ kg/sec$ and  cubic damping constant $c_3= 0.05\ kg\ sec/m^2$.   Let $x_1$ denote the angle of pendulum measured counter clockwise from the outward pointing ray from the center of the Earth and let $x_2$ denote the angular velocity.  The continuous time determistic equations of motion are
\bean
\dot{x}_1&=& x_2
\\
\dot{x}_2&=& lg\sin x_1 -c_1 x_2-c_3 x_2^3 +u
\eean
 The goal is to find a feedback $u=\kappa(x)$ that stabilizes the pendulum to straight up in spite of the noises so we take the   continuous time criterion to be 
\bean
\min_u {1\over 2}\int_0^\infty \|x\|^2+u^2\ dt
\eean

We time discretize this problem using Euler's method with a time step of $0.02$ seconds to  get the discrete time optimal control problem
of minimizing
\bean
\min_u 0.01\sum_{t=0}^\infty \|x\|^2+u^2
\eean
subject to
\bean
x_1^+&=& x_1+0.02x_2\\
x_2^+&=& x_2+0.02\left(  lg\sin x_1 -c_1 x_2-c_3 x_2^3 +u\right)
\eean

But the shape of the earth is not a perfect sphere and its density  is not uniform  so there are fluctuations in the "gravity constant".  We model  these relative fluctuations in the "gravity constant" by $0.1 w_1$ although they are probably much smaller.  There might also be  relative fluctuations in the damping constants modeled by $0.1 w_2$.   
We model these stochastically by two  white noises,
\bean
x_1^+&=& x_1+0.02x_2\\
x_2^+&=& x_2+0.02\left(  lg\sin x_1 -c_1 x_2-c_3 x_2^3 +u\right)\\
&&+0.02\left(0.1w_1 lg\sin x_1  - 0.1w_2(c_1 x_2+c_3x_2^3) \right)
\eean
 
This is  an example about how stochastic models with noise coefficients of order $O(x,u)$ can arise.   If the noise is modeling an uncertain  environment then its coefficients are likely  to be  $O(1)$.  But if it is  the  model  that is  uncetain then noise coefficients are likely to be $O(x,u)$.

The linear coefficients in the dynamics are
\bean
F=\bmt 1&0.02\\0.1740&0.9980 \emt,& G=\bmt 0\\0.02\emt,& \\
Q=\bmt 0.02&0\\0&0.02\emt, & R=0.02,& S=\bmt 0\\0\emt\\
C_1=\bmt 0&0\\0.0174&0\emt,&C_2=\bmt 0&0\\0&-0.0002\emt,& \\
D_1=\bmt 0\\0\emt,& D_2 =\bmt 0\\0\emt &
\eean

The above iteration converges in six steps to the solution of SDARE (\ref{sdare1}, \ref{sdare2}),
\bean
P&=& \bmt 54.9340 & 17.9795\\17.9795&6.0744\emt\\
K&=& \bmt -17.9795  & -6.0744\emt
\eean
The eigenvalues of $F+GK$ are $0.9483$ and $ 0.9282$.

By way of comparison if we delete the noise terms from the problem then the solution to DARE is
\bean
P&=& \bmt 54.8930 & 17.9739\\17.9739&6.0734\emt\\
K&=& \bmt -16.9694  & -5.7253\emt
\eean
and the eigenvalues of $F+GK$ are $0.9510$ and $ 0.9325$.

The dynamics is an odd function of $x,u$ so its
 quadratic and quartic  terms are zero.   The cubic terms  are
\bean
f^{[3]}(x,u)&=& \bmt 0\\ -0.029 x_1^3-0.001 x_2^3 \emt\\
\gamma_1^{[3]}(x,u)&=& \bmt 0\\  -0.0029 x_1^3
\emt\\
\gamma_2^{[3]}(x,u)&=& \bmt 0\\  -0.0001 x_1^3
\emt
\eean
and the quintic terms 
are
\bean
f^{[5]}(x,u)&=& \bmt 0\\ 0.00145 x_1^5 \emt\\
\gamma_1^{[5]}(x,u)&=& \bmt 0\\  0.000145 x_1^5
\emt\\
\gamma_2^{[5]}(x,u)&=& \bmt 0\\  0\emt\\
\eean

Because the Lagrangian is an even function  and the dynamics is an odd function of $x,u$
 we know that
$\pi(x)$ is an even function  of $x$ and $\kappa(x)$ is an odd function of $x$.

We have computed the optimal cost $\pi(x)$ to degree $6$ and the optimal feedback $\kappa(x)$ to degree $5$,
\bean  
\pi(x)&=&27.4670x_1^2+   17.9795x_1x_2    +3.0372x_2^2\\&&
-4.4633x_1^4   -2.7258x_1^3x_2   -0.4995_1^2x_2^2   -0.0796x_1x_2^3  -0.0169x_2^4\\
&&
0.3860x_1^6+    01976x_1^5x_2  +0.0266x_1^4x_2^2   +0.0021_1^3x_2^3 \\
&& -0.0003  x_1^2x_2^4   -0.0001x_1x_2^5    +0.00004 x_2^6
\\
\kappa(x)&=&-17.9795x_1   -6.0744x_2 \\
&& 2.7244x_1^3+    0.9604x_1^2x_2+    0.1913x_1x_2^2+    0.0557x_2^3\\
&&
-0.17347x_1^5+   -0.0359x_1^4x_2+0.0056x_1^3x_2^2\\&&+0.0048x_1^2x_2^3+0.0010x_1x_2^4   -0.0001x_2^5
\eean

\begin{figure}
 \centering
\includegraphics[width=3.5in]{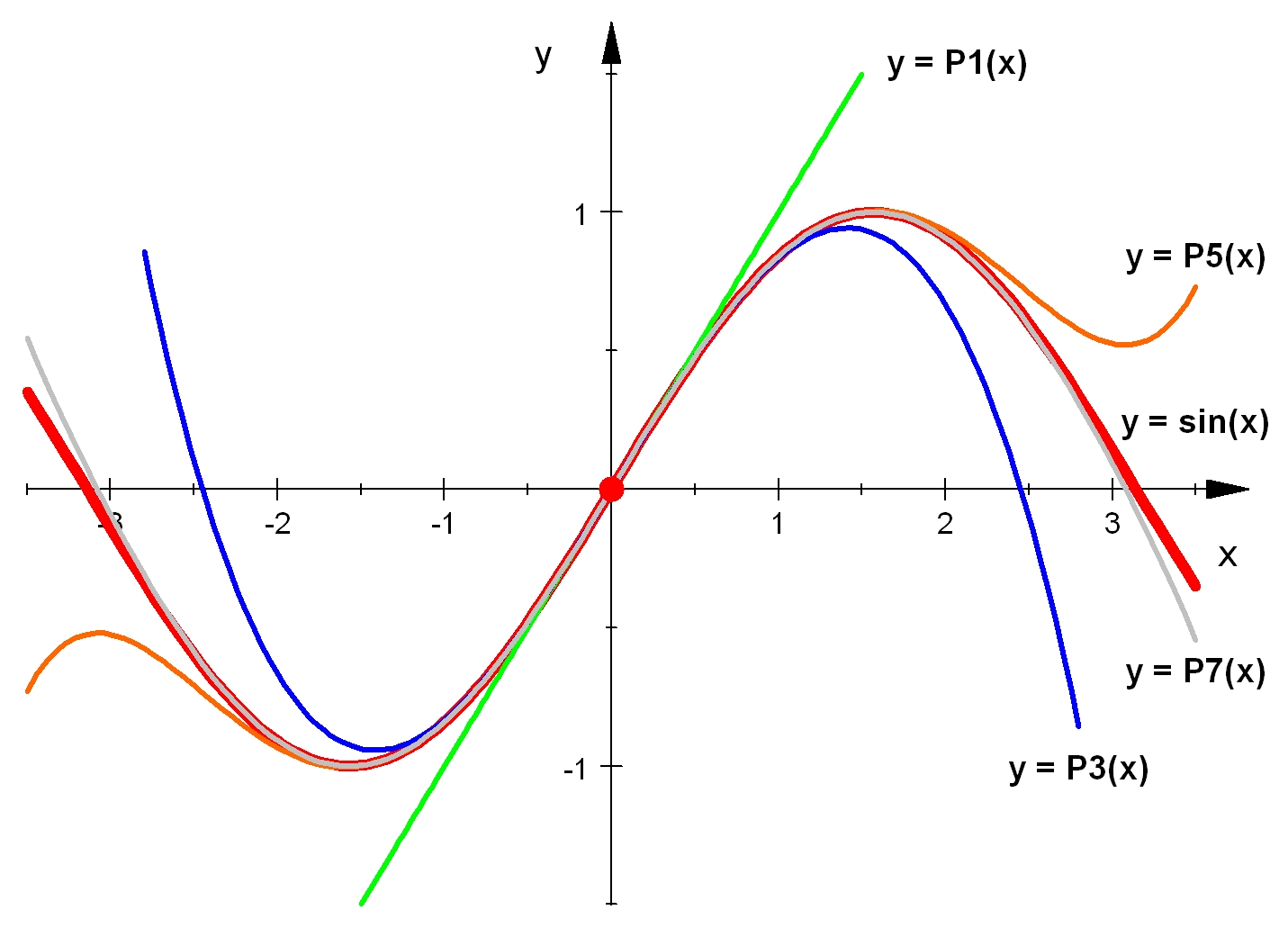}
\caption{Taylor approximations of sin(x)}        
 \end{figure}

In making this computation we are approximating $\sin x_1$ by its Taylor polynomials 
\bean
\sin x_1&=& x_1-{x_1^3\over 6} +{x_1^5 \over 120}+\ldots
\eean
The alternating signs of the odd terms in these polynomials are reflected in the nearly  alternating signs in the Taylor polynomials of the optimal cost $\pi(x)$
and optimal feedback $\kappa(x)$.  If we take a first degree approximation to $\sin x_1$ we are overestimating the gravitational force 
pulling the pendulum from its upright position pointing so $\pi^{[2]}(x)$ overestimates the optimal cost
and the feedback $u=\kappa^{[1]}(x)$ is stronger than it needs to be.  This could be a problem if there is a bound on the magnitude of $u$ that we ignored in the analysis.
If we take a third degree approximation to $\sin x_1$ then $\pi^{[2]}(x)+\pi^{[4]}(x)$ under estimates the optimal cost
and the feedback $u=\kappa^{[1]}(x)+\kappa^{[3]}(x)$ is weaker than it needs to be.
If we take a fifth degree approximation to $\sin x_1$ then $\pi^{[2]}(x)+\pi^{[4]}(x)+\pi^{[6]}(x)$ over estimates the optimal cost but by a smaller margin than 
$\pi^{[2}(x)$.  The feedback     $u=\kappa^{[1]}(x)+\kappa^{[3]}(x)+\kappa^{[5]}(x)$  is stronger than it needs to be
but by a smaller margin than $u=\kappa^{[1]}(x)$.

\section{Finite Horizon Stochastic Nonlinear Optimal Control  Problem} \label{FH}
\setcounter{equation}{0}
Consider the finite horizon stochastic nonlinear optimal control  problem,
\bean
\min_{u(\cdot)}  {\rm E}\left\{ \sum_{t=0}^{T-1} l(t,x,u) +\pi_T(x(T))\right\}
\eean
subject to 
\bean
x^+&=& f(t,x,u) +\sum_{k=1}^r w_k\gamma_k(t,x,u)
\eean
Again we assume that $ f, l,\gamma_k, \pi_T$ are sufficiently smooth.

If they exist and are smooth the optimal cost $\pi(t, x) $ of starting at $x$ at time $t$ and the optimal feedback
$u(t)=\kappa(t,x(t))$ satisfy the Finite Horizon Stochastic Dynamic Programming Equations (FSDPE) (\ref{sdpe1},
\ref{sdpe2})

The quantity to be minimized  is a smooth function of $u$ hence  (\ref{sdpe1}. \ref{sdpe2}) imply
  \bea \nonumber
\pi(t,x)&=& {\rm E} \left\{ \pi\left(t+1,f(t,x,\kappa(t,x))+\sum_{k=1}^rw_k\gamma_k(t,x,\kappa(t,x))\right)\right\}\\
&&+  l(t,x,\kappa(t,x)) \label{sdpe3}\\
\nonumber
0&=& {\rm E} \left\{\frac{\partial \pi}{\partial x}\left(t+1,f(t,x,\kappa(t,x))+\sum_{k=1}^rw_k\gamma_k(t,x,\kappa(t,x)) \right)\right.
\\&& \left.  \times
\left(\frac{\partial f}{\partial u}(t,x,\kappa(t,x)) +\sum_k w_k\frac{\partial  \gamma_k}{\partial u}(t,x,\kappa(t,x))\right)\right\}\nonumber\\
&&+\frac{\partial l}{\partial u}(t,x,\kappa(t,x)) \label{sdpe4}
\eea
Of course the reverse implication is not necessarily true  as the quantity to be minimized could have local minima  or stationary points.

These equations are solved backward in time from the
 final condition
\bea \label{hjbT}
\pi(T,x)&=& \pi_T(x)
\eea

Again we assume that we have the  following Taylor expansions
\bean
f(t,x,u)&=& F(t)x+G(t)u+f^{[2]}(t,x,u)+f^{[3]}(t,x,u)+\ldots\\
l(t,x,u)&=& {1\over 2}\left( x'Q(t)x+2x'S(t)u+u'R(t)u\right)+l^{[3]}(t,x,u)+l^{[4]}(t,x,u)+\ldots\\
\gamma_k(t,x,u)&=& C_k(t)x+D_k(t)u+\gamma_k^{[2]}(t,x,u)+\gamma_{k}^{[3]}(t,x,u)+\ldots\\
\pi_T(x)&=& {1\over 2} x'P_Tx+\pi_T^{[3]}(x)+\pi_T^{[4]}(x)+\ldots\\
\pi(t,x)&=& {1\over 2} x'P(t)x+\pi^{[3]}(t,x)+\pi^{[4]}(t,x)+\ldots\\
\kappa(t,x)&=& K(t)x+\kappa^{[2]}(t,x)+\kappa^{[3]}(t,x)+\ldots
\eean
where $^{[r]}$ indicates terms of homogeneous degree $r$ in $x,u$ with coefficients that are continuous  functions of $t$.
The key assumption is that $\gamma_k(t,0,0)=0$
for then (\ref{sdpe3}, \ref{sdpe4}, \ref{hjbT}) are  amenable to power series methods.

We plug these expansions into the simplified Finite Horizon Stochastic Dynamic Programming Equations(\ref{sdpe3}, \ref{sdpe4}) and collect terms of lowest degree, that is, degree  two in (\ref{sdpe3}, degree one in (\ref{sdpe4}) and degree two in (\ref{hjbT}).
We plug these into SIDPE  which simplifies to
\bea
P(t)&=&  Q(t)+K'(t)S(t) +S(t)K'(t)+K'(t)R(t)K(t) \label{fsdare1} \\
&&+(F(t)+G(t)K(t))'P(t+1)(F(t)+G(t)K(t))  \nonumber\\
&&+\sum_{k=1^r}(C_k(t)+D_k(t)K(t))'P(t+1)(C_k(t)+D_k(t)K(t))  \nonumber 
\\
K(t)&=&-\left(R(t)+G'(t)P(t+1)G(t)+\sum_{k=1}^rD'_k(t)P(t+1)D_k(t)\right)^{-1} \label{fsdare2}
\\ && \times \left(G'(t)P(t+1)F(t)+S'(t)+\sum_{k=1}^rD'_k(t)P(t+1)C_k(t)\right) 
\nonumber
\eea
We call these equations the stochastic discrete time Riccati difference equations (SDRDE).     These difference equations are solved backward in time from the terminal condition
\bean
P(T)&=& P_T
\eean

Then we may proceed to the next degrees, degree three in (\ref{sdpe3}),  and degree two  in (\ref{sdpe4}).
\bea 
&&\pi^{[3]}(t,x)={\rm E} \left\{\pi^{[3]}\left( t+1,z(t,x,w)\right)\right\}\label{sdpe5}\\
&&+x'(F(t)+G(t)K(t))'P(t+1)f^{[2]}(t,x,Kx)\nonumber\\&&
+\sum_k x'(C_k(t)+D_k(t)K(t))'P(t+1)\gamma_k^{[2]}(t,x,Kx)+l^{[3]}(t,x,Kx)\nonumber \\
&&0={\rm E} \left\{\frac{\partial \pi^{[3]}}{\partial x}\left( t,z(t,x,w)\right)
\left(G(t)+\sum_kw_k D_k(t)  \right)  \right\}  \label{sdpe6} 
\\
&&+x'P(t+1)\frac{\partial f^{[2]}}{\partial u}(t,x,K(t)x)+\frac{\partial l^{[3]}}{\partial u}(t,x,K(t)x)  \nonumber \\
&&+(\kappa^{[2]}(t,x))'\left(R(t)+G'(t)P(t+1)G(t)+\sum_kD_k'(t)P(t+1)D_k(t)\right) \nonumber 
\eea
where
\bean
z(t,x,w)&=&F(t)+G(t)K(t))x+\sum_k w_k (C_k(t)+D_k(t)K(t)x
\eean

Notice again the unknown $\kappa^{[2]}(t,x)$ does not appear in the first equation which is linear difference equation for  
$ \pi^{[3]}(t,x)$ running backward in time from the terminal condition,
\bean
\pi^{[3]}(t,x)&=& \pi^{[3]}_T(x)
\eean  
We can solve it and  if $R(t)+G'(t)P(t+1)G(t)+\sum_k D'_k(t)P(t)D_k(t)$ is invertible then we  can solve the second equation for $\kappa^{[2]}(t,x)$. 
 The higher degree terms can be found in a similar fashion.

\end{document}